\documentclass[11pt]{amsart}

\usepackage{amsmath}
\usepackage{amssymb}
\usepackage{mathrsfs}

\newtheorem{theorem}{Theorem}[section]

\newtheorem{lemma}[theorem]{Lemma}

\newtheorem{corollary}[theorem]{Corollary}

\theoremstyle{definition}
\newtheorem{definition}[theorem]{Definition}

\newtheorem{question}[theorem]{Question}

\theoremstyle{remark}

\newcount\skewfactor
\def\mathunderaccent#1#2 {\let\theaccent#1\skewfactor#2
\mathpalette\putaccentunder}
\def\putaccentunder#1#2{\oalign{$#1#2$\crcr\hidewidth
\vbox to.2ex{\hbox{$#1\skew\skewfactor\theaccent{}$}\vss}\hidewidth}}

\def\smallbox#1{\leavevmode\thinspace\hbox{\vrule\vtop{\vbox
   {\hrule\kern1pt\hbox{\vphantom{\tt/}\thinspace{\tt#1}\thinspace}}
   \kern1pt\hrule}\vrule}\thinspace}

\newcommand{\cf}{{\rm cf}}

\usepackage[all]{xy}

\def\qedref#1{$\qed_{\reforiginal{#1}}$}

\title[monochromatic paths]{Infinite monochromatic paths and a theorem of Erd\H{o}s-Hajnal-Rado}
\author{Shimon Garti}
\address{Einstein Institute of Mathematics,
 The Hebrew University of Jerusalem,
 Jerusalem 91904, Israel}
\email{shimon.garty@mail.huji.ac.il}

\author{Menachem Magidor}
\address{Einstein Institute of Mathematics,
 The Hebrew University of Jerusalem,
 Jerusalem 91904, Israel}
\email{menachem@math.huji.ac.il}

\author{Saharon Shelah}
\address{Einstein Institute of Mathematics,
 The Hebrew University of Jerusalem,
 Jerusalem 91904, Israel,
 and Department of Mathematics
 Rutgers University
 New Brunswick, NJ 08854, USA}
\email{shelah@math.huji.ac.il}
\urladdr{http://www.math.rutgers.edu/\char`\~shelah}

\thanks{The first and the third authors are grateful to the support of the ERC, grant no. 338821. This is publication 1165 of the third author}

\begin{document}
\let\labeloriginal\label
\let\reforiginal\ref

\begin{abstract}
We prove that if $\mu$ is a singular cardinal with countable cofinality and $2^\mu=\mu^+$ then $\binom{\mu^+}{\mu}\nrightarrow\binom{\mu^+\ \aleph_2}{\mu\ \mu}$.
\end{abstract}

\maketitle

\newpage

\section{Introduction}

The polarized partition relation $\binom{\alpha}{\beta}\rightarrow \binom{\gamma_0\ \gamma_1}{\delta_0\ \delta_1}$ says that for every coloring $c:\alpha\times\beta\rightarrow 2$ there are $A\subseteq\alpha,B\subseteq\beta$ and $i\in\{0,1\}$ such that ${\rm otp}(A)=\gamma_i,{\rm otp}(B)=\delta_i$ and $c\upharpoonright(A\times B)$ is constantly $i$.
If $\gamma_0=\gamma_1=\gamma$ and $\delta_0=\delta_1=\delta$ then we write $\binom{\alpha}{\beta}\rightarrow\binom{\gamma}{\delta}_2$, in which case we shall say that the relation is \emph{balanced}.

A central case is $\alpha=\mu^+,\beta=\mu$ where $\mu$ is an infinite cardinal.
We shall focus on the subcase in which $\mu$ is a singular cardinal and $2^\mu=\mu^+$.
If $\kappa$ is an infinite cardinal and $2^\kappa=\kappa^+$ then $\binom{\kappa^+}{\kappa}\nrightarrow\binom{\kappa^+}{\kappa}_2$, and this holds at any infinite cardinal including singular cardinals as proved in \cite{MR0202613}.
We say that the \emph{strong polarized relation} fails at the pair $(\kappa,\kappa^+)$ under the local instance of GCH at $\kappa$.

However, in many cases one can force the \emph{almost strong polarized relation} which asserts that $\binom{\mu^+}{\mu}\rightarrow\binom{\alpha}{\mu}_2$ for every $\alpha\in\mu^+$.
If one forces this relation and then collapses $2^\mu$ to $\mu^+$ then one obtains $\binom{\mu^+}{\mu}\rightarrow\binom{\alpha}{\mu}_2$ for every $\alpha\in\mu^+$ with $2^\mu=\mu^+$, see \cite[Claim 3.3]{MR3610266} which is based on an observation of Foreman.
Moreover, if $\mu$ is a limit of measurable cardinals then the positive almost strong relation $\binom{\mu^+}{\mu}\rightarrow\binom{\alpha}{\mu}_2$ for every $\alpha\in\mu^+$ holds in ZFC and hence compatible with $2^\mu=\mu^+$, see \cite{MR1606515}.

These theorems show that we understand the balanced polarized relation quite well at successors of singular cardinals.
A natural question concerns the intermediate \emph{unbalanced} relation $\binom{\mu^+}{\mu}\rightarrow\binom{\mu^+\ \alpha}{\mu\ \mu}$ under the assumption $2^\mu=\mu^+$.
This question has been investigated in \cite{MR0202613}, and the authors proved that if $\mu\geq\cf(\mu)>\aleph_0$ then $2^\mu=\mu^+$ implies the negative relation $\binom{\mu^+}{\mu}\nrightarrow\binom{\mu^+\ \omega}{\mu\ \mu}$.
Namely, the unbalanced relation behaves much similarly to the strong polarized relation, yielding a negative statement under the assumption $2^\mu=\mu^+$ already when $\alpha=\omega$ in the second color.
We indicate that the authors of \cite{MR0202613} assume GCH for this result, but only the local instance at $\mu$ is needed for the proof.

Surprisingly, if $\cf(\mu)=\omega$ then $\binom{\mu^+}{\mu}\rightarrow\binom{\mu^+\ \omega}{\mu\ \mu}$ even if $2^\mu=\mu^+$.
This is an extremely rare situation, in which singular cardinals with countable cofinality demonstrate a better combinatorial relation than cardinals with uncountable cofinality.
Nevertheless, Erd\H{o}s, Hajnal and Rado proved that if the GCH holds then $\binom{\mu^+}{\mu}\nrightarrow\binom{\mu^+\ \omega_2}{\mu\ \mu}$ whenever $\mu>\cf(\mu)=\omega$.
Unlike the case of $\cf(\mu)>\omega$, here $2^\mu=\mu^+$ is insufficient for their argument.
The GCH assumption can be relaxed, but a crucial assumption used within the proof is that $2^\omega=\omega_1$.
Under this assumption $\omega_2=(2^\omega)^+$ and hence the celebrated Erd\H{o}s-Rado theorem holds at $\omega_2$, namely $\omega_2\rightarrow(\omega_1)^2_\omega$.
This theorem is essential for proving the negative polarized relation $\binom{\mu^+}{\mu}\nrightarrow\binom{\mu^+\ \omega_2}{\mu\ \mu}$.

The employment of the Erd\H{o}s-Rado theorem invites for a natural question.
It is simple to see that if $\omega<\kappa\leq 2^\omega$ then $\kappa\nrightarrow(\omega)^2_\omega$ and acutally $\kappa\nrightarrow(3)^2_\omega$.
One may wonder, therefore, whether $\binom{\mu^+}{\mu}\nrightarrow\binom{\mu^+\ \kappa}{\mu\ \mu}$ where $\kappa\leq 2^\omega$ and $2^\mu=\mu^+$.
This question can be phrased at two levels.
In the absolute level one may ask about $\kappa=\omega_1$, and in the non-absolute level at $\kappa\leq 2^\omega$ when $2^\omega>\omega_1$ is forced.

\begin{question}
\label{qehr} Assume that $2^\mu=\mu^+$ and $\mu>\cf(\mu)=\omega$.
\begin{enumerate}
\item [$(\aleph)$] Is it consistent that $\binom{\mu^+}{\mu}\rightarrow\binom{\mu^+\ \omega_1}{\mu\ \mu}$?
\item [$(\beth)$] Is it consistent that $\binom{\mu^+}{\mu}\rightarrow\binom{\mu^+\ \kappa}{\mu\ \mu}$ for some $\omega<\kappa\leq 2^\omega$?
\end{enumerate}
\end{question}

Regarding the first part of the question we indicate that Jones proved in \cite{MR2367118} that $\binom{\mu^+}{\mu}\rightarrow\binom{\mu^+\ \alpha}{\mu\ \mu}$ for every countable ordinal $\alpha$, thus reaching well-nigh to $\omega_1$.
Our main result here is a sharp negative answer to the second part.
Namely, if $\mu>\cf(\mu)=\omega$ then $\binom{\mu^+}{\mu}\nrightarrow\binom{\mu^+\ \omega_2}{\mu\ \mu}$ under the assumption $2^\mu=\mu^+$, no matter how large is $2^\omega$.
This includes, in particular, cases in which $2^\mu=\mu^+$ but $\mu$ is not strong limit.

In order to prove our result we replace the Erd\H{o}s-Rado theorem by a statement about monochromatic paths in complete graphs.
We shall use a theorem of Todor\v{c}evi\'c from \cite{MR793235}, see also \cite{MR1086126}.
It says that $\omega_2\rightarrow_{\rm ipath}(\omega)^2_\omega$ even if $2^\omega\geq\omega_2$.
On the other hand, $\omega_1\nrightarrow_{\rm ipath}(\omega)^2_\omega$ and even $\omega_1\nrightarrow_{\rm ipath}(\omega)^2_{\omega,<\omega}$ as shown by Todor\v{c}evi\'c in \cite{MR908147}.
Hence concerning the polarized relation under scrutiny, we do not resolve the focal case of $\omega_1$.

Along the refereeing process we learned that the above mentioned path relations were first proved by Todor\v{c}evi\'c.
Our original manuscript contained similar statements with our proofs.
Since the proofs are quite different (and some of the statements are not identical) we include a discussion on path relations in the current version.
The paper is organized in such a way that the new result (concerning the polarized relation) appears in the first section, and path relations are discussed in the second, so the reader may skip that part.
We are deeply indebted to the referees for pointing to the literature concerning path relations, and for many other helpful suggestions.

Our notation is mostly standard.
We follow \cite{MR795592} with respect to arrows notation.
Our set theoretical notation is coherent, in general, with \cite{MR1940513}, but we adopt the Jerusalem notation in forcing, so $p\leq_{\mathbb{P}}q$ reads $p$ is weaker than $q$.
An excellent background concerning the basics of the polarized relation can be found in \cite{MR3075383}.

\newpage

\section{Polarized relations}

Assume that $\theta\leq\kappa$.
We shall say that $\kappa\rightarrow_{\rm ipath}(\omega)^2_\theta$ iff for every $c:[\kappa]^2\rightarrow\theta$ there exist an increasing sequence $(\beta_n:n\in\omega)$ of ordinals of $\kappa$ and a color $\gamma\in\theta$ so that $c(\beta_n,\beta_{n+1})=\gamma$ for every $n\in\omega$.
We shall say that $\kappa\rightarrow_{\rm ipath}(\omega)^2_{\theta,<\omega}$ iff for every $c:[\kappa]^2\rightarrow\theta$ there exists an increasing sequence $(\beta_n:n\in\omega)$ of ordinals of $\kappa$ such that $|\{c(\beta_n,\beta_{n+1}):n\in\omega\}|<\aleph_0$.
These relations are weakenings of the classical Erd\H{o}s-Rado theorem.
In particular, if $\theta=\omega$ and $\kappa=(2^\omega)^+$ then $\kappa\rightarrow_{\rm ipath}(\omega)^2_\theta$.

Several mathematicians considered these relations, and it seems that Silver, \cite{MR274278}, was the first one to use them.
Todor\v{c}evi\'c proved the following two statements, among many other, in \cite{MR793235} and \cite{MR908147}.
We refer also to \cite[Corollary 14]{MR1086126} in this context.

\begin{theorem}
\label{thmtod} [Todor\v{c}evi\'c] Path relations.
\begin{enumerate}
\item [$(\aleph)$] $\omega_2\rightarrow_{\rm ipath}(\omega)^2_\omega$.
\item [$(\beth)$] $\omega_1\nrightarrow_{\rm ipath}(\omega)^2_{\omega,<\omega}$.
\end{enumerate}
\end{theorem}

\hfill \qedref{thmtod}

Our purpose in this section is to derive negative polarized relations at singular cardinals with countable cofinality from instances of the ordinary path relation.
We shall need the following:

\begin{definition}
\label{defalternative} Polarized relations with alternatives. \newline
We say that $\begin{pmatrix} \alpha \\ \beta \end{pmatrix} \rightarrow
\begin{pmatrix}
\begin{matrix} \gamma_0 \\ \delta_0 \end{matrix}
\vee
\begin{matrix} \varepsilon_0\ \\ \zeta_0\ \end{matrix}
\begin{matrix} \gamma_1 \\ \delta_1 \end{matrix}
\vee
\begin{matrix} \varepsilon_1 \\ \zeta_1 \end{matrix}
\end{pmatrix}$ iff for every coloring $c:\alpha\times\beta\rightarrow 2$ there are $i\in\{0,1\}, A\subseteq\alpha$ and $B\subseteq\beta$ such that $c''(A\times B)=\{i\}$ and either ${\rm otp}(A)=\gamma_i, {\rm otp}(B)=\delta_i$ \emph{or} ${\rm otp}(A)=\varepsilon_i, {\rm otp}(B)=\zeta_i$.
\end{definition}

It follows from the definition that if some relation holds with alternatives then it holds upon omitting one of the alternatives (or more).
Of course, one may suggest an alternative only in one of the colors, as done in the following theorem which is the main result of this section.

\begin{theorem}
\label{thmmt} Assume that:
\begin{enumerate}
\item [$(\aleph)$] $\mu>\cf(\mu)=\omega$.
\item [$(\beth)$] $2^\mu=\mu^+$.
\item [$(\gimel)$] $\kappa\rightarrow_{\rm ipath}(\omega)^2_\omega$, or even $\kappa\rightarrow_{\rm ipath}(\omega)^2_{\omega,<\omega}$.
\end{enumerate}
Then $\begin{pmatrix} \mu^+ \\ \mu^+ \end{pmatrix} \nrightarrow \begin{pmatrix} \begin{matrix} \mu^+\ \\ \mu\ \end{matrix} \begin{matrix} \kappa \\ \mu \end{matrix} \vee \begin{matrix} 1 \\ \mu^+ \end{matrix} \end{pmatrix}$.
\end{theorem}

\par\noindent\emph{Proof}. \newline
Let $\{B_\alpha:\alpha\in\mu^+\}$ enumerate the elements of $[\mu^+]^\mu$.
For every $0<\alpha\in\mu^+$ we reenumerate the family $\{B_\beta:\beta\in\alpha\}$ by $\{B_{\alpha\varepsilon}:\varepsilon\in\mu\}$ (possibly with repetitions).
For every $0<\alpha\in\mu^+$ we also reenumerate the ordinals of $\alpha$ by $\{\alpha_\eta:\eta\in\mu\}$, so our enumerations are of order type $\mu$.

Now for each $\alpha\in\mu^+$ and every $\varepsilon\in\mu$ we pick up an ordinal $\gamma(\alpha,\varepsilon)\in\mu^+$ such that $\gamma(\alpha,\varepsilon)\in B_{\alpha\varepsilon} - \{\gamma(\alpha_\eta,\zeta):\eta<\varepsilon,\zeta\leq\varepsilon\}$.
The choice is possible since $\varepsilon\in\mu$ and $|B_{\alpha\varepsilon}|=\mu$.
Based on these ordinals we define a coloring $c:\mu^+\times\mu^+\rightarrow 2$ as follows:
$$
c(\alpha,\beta)=1 \Leftrightarrow \exists\varepsilon\in\mu, \beta=\gamma(\alpha,\varepsilon).
$$
In order to show that $c$ exemplifies the alleged negative relation we show, first, that there is no $0$-monochromatic product of size $\mu^+\times\mu$ under $c$.
Assume, therefore, that $A,B\subseteq\mu^+, |A|=\mu^+$ and $|B|=\mu$.
Pick up a sufficiently large ordinal $\alpha\in A$ so that $B\in\{B_\beta:\beta<\alpha\}$ and let $\varepsilon\in\mu$ be an ordinal for which $B = B_{\alpha\varepsilon}$.
Recall that $\gamma(\alpha,\varepsilon)\in B_{\alpha\varepsilon}=B$, so letting $\beta=\gamma(\alpha,\varepsilon)$ we see that $c(\alpha,\beta)=1$.
It follows that $c\upharpoonright(A\times B)$ is not $0$-monochromatic.

Secondly, we show that there is no $1$-monochromatic product of size $1\times\mu^+$ under $c$.
For this end, suppose that $\alpha\in\mu^+, H\subseteq\mu^+$ and $c\upharpoonright(\{\alpha\}\times H)$ is $1$-monochromatic. We must show that the cardinality of $H$ is less than $\mu^+$.
Notice that $H=\{\beta\in\mu^+:c(\alpha,\beta)=1\}\subseteq \{\gamma(\alpha,\varepsilon):\varepsilon\in\mu\}$, hence $|H|\leq|\{\gamma(\alpha,\varepsilon):\varepsilon\in\mu\}|\leq\mu$ as required.

Finally, let us show that there is no $1$-monochromatic product of size $\kappa\times\mu$ under $c$.
Observe that for every $\alpha\in\mu^+$ we have:
$$
\eta<\varepsilon<\mu\wedge\gamma(\alpha_\eta,\zeta)=\gamma(\alpha,\varepsilon) \Rightarrow \varepsilon<\zeta.
$$
Indeed, $\alpha_\eta<\alpha$ and hence if $\zeta\leq\varepsilon$ then in the choice of $\gamma(\alpha,\varepsilon)$ we make sure that $\gamma(\alpha,\varepsilon)\neq\gamma(\alpha_\eta,\zeta)$.
So assuming that $\gamma(\alpha_\eta,\zeta)=\gamma(\alpha,\varepsilon)$ we conclude that $\varepsilon<\zeta$ and hence the above statement holds.
We denote this observation by $(\ast)_\alpha$ for every $\alpha\in\mu^+$.

Assume that $S\subseteq\mu^+, |S|<\mu$ and $S$ is infinite.
We may concentrate on the first $\omega$ elements of $S$, so assume that $S=\{\gamma_n:n\in\omega\}$.
Assume, further, that $\sup\{\eta\in\mu:\exists i\in\omega, \gamma_i=(\gamma_{i+1})_\eta\}<\mu$.
Define $T(S) = \{\beta\in\mu^+:\forall\alpha\in S, c(\alpha,\beta)=1\}$.
We claim that $|T(S)|<\mu$ as well.

For proving this claim fix an ordinal $\rho<\mu$ such that $\sup\{\eta\in\mu:\exists i\in\omega, \gamma_i=(\gamma_{i+1})_\eta\}<\rho$.
Choose any ordinal $\gamma\in T(S)$.
By the very definition of $T(S)$ we see that $c(\alpha,\gamma)=1$ for every $\alpha\in S$.
By the definition of the coloring $c$, for every $\alpha\in S$ there is a unique ordinal $\varepsilon(\alpha)\in\mu$ such that $\gamma=\gamma(\alpha,\varepsilon(\alpha))$.
We claim that for some $\alpha\in S$ we have $\varepsilon(\alpha)<\rho$.

Assume towards contradiction that the claim fails.
Choose a pair of consecutive ordinals $\alpha,\alpha'\in S$ such that $\alpha<\alpha'$.
Fix an ordinal $\eta\in\mu$ for which $\alpha'_\eta=\alpha$.
Notice that $\eta<\rho$, by the definition of $\rho$.
Now the assumption toward contradiction implies that $\eta<\varepsilon(\alpha')$.
Likewise, $\gamma=\gamma(\alpha,\varepsilon(\alpha)) = \gamma(\alpha',\varepsilon(\alpha'))$.
Applying $(\ast)_{\alpha'}$ we conclude that $\varepsilon(\alpha')<\varepsilon(\alpha)$.
Since the set $S$ is infinite, if we choose an increasing sequence $(\alpha_n:n\in\omega)$ of elements of $S$ we produce an infinite decreasing sequence of ordinals $(\varepsilon(\alpha_n):n\in\omega)$.
This absurd confirms the claim.

Hence for every $\gamma\in T(S)$ we fix an ordinal $\alpha\in S$ such that $\gamma=\gamma(\alpha,\varepsilon(\alpha))$ and $\varepsilon(\alpha)<\rho$.
It follows that $T(S)\subseteq\{\gamma(\alpha,\varepsilon): \varepsilon<\rho,\alpha\in S\}$.
Consequently, $|T(S)|\leq|\rho|\cdot|S|<\mu$.

Suppose that $S\subseteq\mu^+, |S|=\kappa$ and $c\upharpoonright(S\times T)$ is $1$-monochromatic.
This means that $T\subseteq T(S)$, so suffice it to show that $|T(S)|<\mu$.
If we could only show that $\sup\{\eta\in\mu: \exists i\in\omega, \gamma_i=(\gamma_{i+1})_\eta\}<\mu$ then using the fact that $|S|=\kappa<\mu$ we will be done.
In order to prove this bound, let $(\mu_n:n\in\omega)$ be an increasing sequence of regular cardinals such that $\mu=\bigcup_{n\in\omega}\mu_n$.
Define a coloring $d:[S]^2\rightarrow\omega$ as follows:
$$
d(\gamma,\delta)=n \Leftrightarrow \gamma<\delta\wedge\exists\eta<\mu_n, \gamma=\delta_\eta.
$$
From the assumption $\kappa\rightarrow_{\rm path}(\omega)^2_\omega$ we infer that there is an increasing sequence $(\gamma_n:n\in\omega)$ of ordinals in $S$ and a color $m\in\omega$ such that $d(\gamma_n,\gamma_{n+1})=m$ for every $n\in\omega$.
Denote the set $\{\gamma_n:n\in\omega\}$ by $S_0$.
Since $S_0\subseteq S$ we see that $T(S_0)\supseteq T(S)$, so it is sufficient to prove that $|T(S_0)|<\mu$.
By the definition of $d$ we see that $\sup\{\eta\in\mu: \exists i\in\omega, \gamma_i=(\gamma_{i+1})_\eta\}\leq\mu_m<\mu$ and even if we assume only $\kappa\rightarrow_{\rm path}(\omega)^2_{\omega,<\omega}$ then $\sup\{\eta\in\mu: \exists i\in\omega, \gamma_i=(\gamma_{i+1})_\eta\}\leq\mu_m<\mu$ for some $m\in\omega$.
Hence $|T(S_0)|<\mu$ and the proof is accomplished.

\hfill \qedref{thmmt}

We can address now the second part of Question \ref{qehr} by eliminating the assumption $2^\omega=\omega_1$ from the proof of $\binom{\mu^+}{\mu}\nrightarrow\binom{\mu^+\ \omega_2}{\mu\ \mu}$ which appears in \cite{MR0202613}.

\begin{corollary}
\label{corehr} Suppose that $\mu>\cf(\mu)=\omega$ and $2^\mu=\mu^+$. \newline
Then $\binom{\mu^+}{\mu}\nrightarrow\binom{\mu^+\ \kappa}{\mu\ \mu}$ and even $\binom{\mu^+}{\mu^+}\nrightarrow\binom{\mu^+\ \kappa}{\mu\ \mu}$ whenever $\omega_1<\kappa$.
\end{corollary}

\hfill \qedref{corehr}

The problem of $\binom{\mu^+}{\mu}\rightarrow\binom{\mu^+\ \omega_1}{\mu\ \mu}$ at singular cardinals with countable cofinality under the assumption $2^\mu=\mu^+$ remains open.
We believe that a positive relation is consistent.
Moreover, in the light of \cite{MR1606515} we even raise the possibility that it holds in ZFC under sufficiently strong assumptions of large cardinals:

\begin{question}
\label{q586} Suppose that $\mu$ is an $\omega$-limit of measurable cardinals. \newline
Is it provable that $\binom{\mu^+}{\mu}\rightarrow\binom{\mu^+\ \omega_1}{\mu\ \mu}$?
\end{question}

\newpage

\section{Path relations the hard way}

In this section we consider path relations at $\omega_2$ and $\omega_1$.
There is a slight difference between our concepts and the path relations of Scheepers and Todor\v{c}evi\'c mentioned in the previous section, as will be explicated anon.
We emphasize that for the main result of the paper concerning polarized relations one can use the statements of Todor\v{c}evi\'c.

Path relations were considered in \cite{MR274278} and in \cite{MR450075}, as well as in \cite{MR0299537}.
The latter has been continued in \cite{MR1050558}, \cite{MR896028} and \cite{MR1050560}.
Due to the terminology of \cite{MR1050558} we distinguish two notions of paths.
If the elements of the path are ordinals and they appear in increasing order then we shall say \emph{an increasing path}.
If the elements of the path are ordinals with no requirement of being increasing then we shall say \emph{a simple path}.
We use the notation $\kappa\rightarrow_{\rm ipath}(\omega)^2_\theta$ in the first case, and $\kappa\rightarrow_{\rm path}(\omega)^2_\theta$ in the second.
Scheepers and Todor\v{c}evi\'c dealt with increasing paths, while the paths in this section are simple.

\begin{definition}
\label{defpath} Simple path relations. \newline
Let $\theta$ and $\kappa$ be cardinals.
\begin{enumerate}
\item [$(\aleph)$] The relation $\kappa\rightarrow_{\rm path}(\omega)^2_\theta$ holds iff for every $c:[\kappa]^2\rightarrow\theta$ there exist a sequence $(\beta_n:n\in\omega)$ of distinct ordinals of $\kappa$ and a color $\gamma\in\theta$ so that $c(\beta_n,\beta_{n+1})=\gamma$ for every $n\in\omega$.
\item [$(\beth)$] The relation $\kappa\rightarrow_{\rm path}(\omega)^2_{\theta,<\omega}$ holds iff for every $c:[\kappa]^2\rightarrow\theta$ there exists a sequence $(\beta_n:n\in\omega)$ of distinct ordinals of $\kappa$ such that $|\{c(\beta_n,\beta_{n+1}):n\in\omega\}|<\aleph_0$.
\end{enumerate}
\end{definition}

Here is our first theorem:

\begin{theorem}
\label{thmomega2}:
\begin{enumerate}
\item [$(\aleph)$] $\omega_2\rightarrow_{\rm path}(\omega)^2_\omega$.
\item [$(\beth)$] If $\kappa\leq 2^\omega$ then there is a coloring of the ordered pairs of $\kappa$ with no monochromatic infinite path.
\item [$(\gimel)$] For every coloring of the ordered pairs of $\omega_2$ there is an infinite path which assumes at most two colors.
\end{enumerate}
\end{theorem}

\par\noindent\emph{Proof}. \newline
Beginning with the first part, let $c:[\omega_2]^2\rightarrow\omega$ be a coloring, and let $\chi$ be a sufficiently large regular cardinal.
Choose an elementary submodel $M\prec(\mathcal{H}(\chi),\in)$ of size $\aleph_1$ so that $\omega_1\subseteq M$ and $\omega_1,c\in M$.
Let $\delta=\omega_2\cap M$ be the characteristic ordinal, so $\delta\in\omega_2$ and we choose $M$ in such a way that $\cf(\delta)=\omega_1$.

For every $n\in\omega$ let $B_n=\{\beta\in\delta:c(\beta,\delta)=n\}$.
Notice that $\delta=\bigcup_{n\in\omega}B_n$, so there exists $n_0\in\omega$ such that $B_{n_0}$ is unbounded in $\delta$.
By induction on $i\in\omega$ we choose an ordinal $\beta_i\in\delta$ such that the following requirements hold for every $i\in\omega$:
\begin{enumerate}
\item [$(a)$] $\beta_{2i}\in B_{n_0}$.
\item [$(b)$] $\beta_{2i+1}>\beta_{2i}$.
\item [$(c)$] $\beta_{2i}>\beta_{2j+1}$ for every $j<i$.
\item [$(d)$] $c(\beta_{2j},\beta_{2i+1})=n_0$ for every $j\leq i$.
\end{enumerate}
The choice is possible since $B_{n_0}$ is unbounded in $\delta$ and since $\cf(\delta)=\omega_1$.
For part $(d)$ we use elementarity.
Now the sequence $(\beta_{2i},\beta_{2i+3}:i\in\omega)$ forms a monochromatic path where consecutive ordinals are colored by $n_0$, thus part $(\aleph)$ has been proved.

For part $(\beth)$ fix an uncountable cardinal $\kappa\leq 2^\omega$.
Let $\mathscr{T}$ be the full binary tree of height $\omega$, and let $(b_\alpha:\alpha\in\kappa)$ enumerate $\kappa$-many distinct $\omega$-branches of $\mathscr{T}$.
We shall define a coloring $c$ over the ordered pairs of $\kappa$ using $\omega+\omega$ colors.
Given two distinct ordinals $\alpha,\beta\in\kappa$ let $m=m(\alpha,\beta)\in\omega$ be the departure level of $b_\alpha$ and $b_\beta$.
Namely, $b_\alpha\upharpoonright m = b_\beta\upharpoonright m$ but $b_\alpha(m)\neq b_\beta(m)$.
Now let $c(\alpha,\beta)=m_0$ iff $b_\alpha(m)=0\wedge b_\beta(m)=1$ and let $c(\alpha,\beta)=m_1$ iff $b_\alpha(m)=1\wedge b_\beta(m)=0$.
Notice that if $c(\alpha,\beta)=m_\ell$ then $c(\beta,\alpha)=m_{1-\ell}$, so the order is crucial here.

It follows that there is no monochromatic infinite path under $c$, and actually no $3$-path which is monochromatic.
Indeed, if $\alpha,\beta,\gamma\in\kappa$ and $c(\alpha,\beta)=c(\beta,\gamma)$ then $b_\alpha\upharpoonright m = b_\beta\upharpoonright m = b_\gamma\upharpoonright m$ for some $m\in\omega$.
But then necessarily $\{c(\alpha,\beta),c(\beta,\gamma)\} = \{m_0,m_1\}$.

We move to part $(\gimel)$ which says that one can limit the above negative examples to two colors only along the path.
To see this, choose $M\prec(\mathcal{H}(\chi),\in)$ as in the proof of part $(\aleph)$, and let $\delta,B_{n_0}$ be as there.

For every $m\in\omega$ let $B_{n_0}^m=\{\beta\in B_{n_0}:c(\delta,\beta)=m\}$ and pick up $m_0\in\omega$ for which $B^{m_0}_{n_0}$ is unbounded in $\delta$.
By induction on $i\in\omega$ we choose $\beta_i$ as in the first part of the proof, but for each $i\in\omega$ we add the requirement that $c(\beta_{2i+1},\beta_{2j})=m_0$ for every $j\leq i$.
Elementarity guarantees that this is possible, recalling the definition of the set $B^{m_0}_{n_0}$.
As before, the path will be $(\beta_i:i\in\omega)$ and one can verify that $\{c(\beta_i,\beta_{i+1}), c(\beta_{i+1},\beta_i):i\in\omega\} = \{m_0,n_0\}$ as required.

\hfill \qedref{thmomega2}

If $2^\omega=\omega_1$ then the above statements are immediate, so the main point of the above theorem is that it holds in ZFC.
In particular, $\omega_2\rightarrow_{\rm path}(\omega)^2_\omega$ holds even if $2^\omega\geq\omega_2$.
The following theorem and corollary show that such a relation is impossible when $\omega_1$ is deemed.
Moreover, even weak homogeneity is excluded.
We emphasize that this result follows from \cite[6.8]{MR908147}, as pointed out by one of the referees. Our method of proof is to employ a forcing argument, and then to claim that it holds in ZFC by absoluteness.

\begin{theorem}
\label{thmomega1} (ZFC): \newline
There exists a coloring $c:[\omega_1]^2\rightarrow\omega$ such that for every finite subset $\{\alpha_0,\ldots,\alpha_m\}\subseteq\omega_1$ and every $v\subseteq\{0,\ldots,m\}$ such that $\langle\alpha_\ell:\ell\in v\rangle$ is an increasing sequence, it is true that $|v|<\max\{c(\alpha_\ell,\alpha_{\ell+1}):\ell<m\}$.
\end{theorem}

\par\noindent\emph{Proof}. \newline
We define a forcing notion $\mathbb{Q}$.
A condition $p\in\mathbb{Q}$ is a pair $(u,f)=(u^p,f^p)$ so that $u\in[\omega_1]^{<\omega}, f:[u]^2\rightarrow\omega$ and if $v\subseteq\{0,\ldots,m\}, \{\alpha_0,\ldots,\alpha_m\}\subseteq u$ and $(\alpha_\ell:\ell\in v)$ is increasing then $|v|<\max\{f(\alpha_\ell,\alpha_{\ell+1}):\ell<m\}$.
For the forcing order, if $p,q\in\mathbb{Q}$ then $p\leq_{\mathbb{Q}}q$ iff $p\subseteq q$, that is $u^p\subseteq u^q$ and $f^q\upharpoonright[u^p]^2=f^p$.
Intuitively, conditions in $\mathbb{Q}$ are finite approximations to the function which we try to force.

For every $\alpha\in\omega_1$ let $\mathcal{D}_\alpha = \{p\in\mathbb{Q}: \alpha\in u^p\}$.
Let us show that $\mathcal{D}_\alpha$ is dense for every $\alpha\in\omega_1$.
Suppose that $\alpha\in\omega_1$ and $p\notin\mathcal{D}_\alpha$, namely $\alpha\notin u^p$.
Let $u^q=u^p\cup\{\alpha\}, f^q\upharpoonright[u^p]^2=f^p$.
If $\beta\in u^p$ then let $f^q(\alpha,\beta) = |u^p|+1+|u^p\cap\beta|$.
One can verify that $q=(u^q,f^q)$ satisfies $p\leq q\in\mathcal{D}_\alpha$ so $\mathcal{D}_\alpha$ is dense.

In Lemma \ref{lemccc} below we shall prove that $\mathbb{Q}$ is $ccc$, so forcing with $\mathbb{Q}$ preserves cardinals.
Let $G\subseteq\mathbb{Q}$ be $V$-generic.
Define $c=\bigcup\{f^p:p\in G\}$.
By the density of each $\mathcal{D}_\alpha$ we see that ${\rm dom}(c)=[\omega_1]^2$.
By the directness of $G$ we see that $c$ is a function.
By the definition of the conditions we see that $c$ exemplifies the required statement.

The coloring $c$ has been forced, but we argue that such a coloring already exists in ZFC.
To see this, notice that the existence of our coloring is expressible as an existence statement of a model of some formula $\psi\in L_{\omega_1\omega}(Q)$ where $Q$ is the quantifier of \emph{there exist uncountably many}.
Indeed, the statement of the theorem asserts that there is a coloring over an uncountable domain.
This can be expressed using the quantifier $Q$.
Now for every $m\in\omega$ we can express the property stated in the theorem by a first order formula.
Using the infinitary logic $L_{\omega_1\omega}$ we can form a conjuction of these statements for every $m\in\omega$, so there is a formula $\psi\in L_{\omega_1\omega}(Q)$ as required.
By \cite{MR0344115} we conclude that such a coloring exists in the ground model, so we are done.

\hfill \qedref{thmomega1}

Recall that a forcing notion $\mathbb{Q}$ is $ccc$ iff $|\mathcal{A}|\leq\aleph_0$ whenever $\mathcal{A}$ is an antichain of $\mathbb{Q}$.
In order to accomplish the proof we need the following:

\begin{lemma}
\label{lemccc} The forcing notion $\mathbb{Q}$ is $ccc$.
\end{lemma}

\par\noindent\emph{Proof}. \newline
We commence with a general claim about projecting conditions in $\mathbb{Q}$ to a countable ordinal.
Suppose that $0<\delta<\omega_1$ and $\delta$ is a limit ordinal.
Suppose that $q\in\mathbb{Q}$ and let $p=q\upharpoonright\delta$, that is $u^p=u^q\cap\delta$ and $f^p=f^q\upharpoonright[u^p]^2$.
Notice that $p=(u^p,f^p)\in\mathbb{Q}$ and $p\leq_{\mathbb{Q}}q$.

Choose an increasing function $h:u^q\rightarrow\delta$ such that $h\upharpoonright u^p$ is the identity function.
Notice that $h''u^q$ is bounded in $\delta$ since $\delta$ is a limit ordinal.
We shall define a condition $p_{[\delta]}$ as follows.
First, let $u^{p_{[\delta]}}$ be $h''q$.
Second, if $\varepsilon,\zeta\in u^{p_{[\delta]}}$ and $\alpha,\beta$ satisfy $\varepsilon=h(\alpha),\zeta=h(\beta)$ then let $f^{p_{[\delta]}}(\varepsilon,\zeta)=f^q(\alpha,\beta)$.
We indicate that $p_{[\delta]}=(u^{p_{[\delta]}},f^{p_{[\delta]}})\in\mathbb{Q}$ and $p\leq p_{[\delta]}$.

We argue that $p_{[\delta]}\parallel q$.
Let us justify this statement and then explain how to derive the chain condition from it.
Our purpose is to define a condition $r$ so that $q,p_{[\delta]}\leq r$.
Let $u^r=u^{p_{[\delta]}}\cup u^q$.
Let $f^r\upharpoonright[u^{p_{[\delta]}}]^2=f^{p_{[\delta]}}$ and $f^r\upharpoonright[u^q]^2=f^q$, upon noticing that on the common part of $u^p$ we assign the same values and on the rest we have disjoint sets so the information is not contradictory.
It remains to define $f^r$ over mixed pairs, so assume that $\alpha\in u^{p_{[\delta]}}-u^p$ and $\beta\in u^q-u^p$.
We let $f^r(\alpha,\beta)=|u^r|+2+|u^r|\cdot(|u^r\cap\alpha|+|u^r\cap\beta|)$.
Finally, let $r=(u^r,f^r)$.

It is clear that $p_{[\delta]},q\leq r$ once we show that $r\in\mathbb{Q}$, so we must prove that $r$ satisfies the defining property of our conditions.
Suppose, therefore, that $v\subseteq\{0,\ldots,m\}, \{\alpha_\ell:\ell\leq m\}\subseteq u^r$ and $\langle\alpha_\ell:\ell\in v\rangle$ is strictly increasing.
We argue that $|v|<\max\{f^r(\alpha_\ell,\alpha_{\ell+1}):\ell<m\}$.
In order to prove this, we consider three cases. \newline

\par\noindent\emph{Case 1}: For some $\ell<m, \{\alpha_\ell,\alpha_{\ell+1}\}\nsubseteq u^{p_{[\delta]}}\wedge \{\alpha_\ell,\alpha_{\ell+1}\}\nsubseteq u^q$.

In this case we see that
$|v|\leq|u^r|<f^r(\alpha_\ell,\alpha_{\ell+1})\leq \max\{f^r(\alpha_\ell,\alpha_{\ell+1}):\ell<m\}$ by the definition of $f^r$ over mixed pairs. \newline

We are left with all the cases in which for every $\ell<m$ either
$\{\alpha_\ell,\alpha_{\ell+1}\}\subseteq u^{p_{[\delta]}}$ or $\{\alpha_\ell,\alpha_{\ell+1}\}\subseteq u^q$.
Hence from now on we assume that Case 1 fails.
Let us remark that it is fairly possible that for some $k,\ell<m$ we will have $\{\alpha_k,\alpha_{k+1}\}\subseteq u^{p_{[\delta]}} \wedge \{\alpha_k,\alpha_{k+1}\}\nsubseteq u^q$ while $\{\alpha_\ell,\alpha_{\ell+1}\}\subseteq u^q \wedge \{\alpha_\ell,\alpha_{\ell+1}\}\nsubseteq u^{p_{[\delta]}}$.
For example, if $\alpha_0\in u^p, \alpha_1\in u^q-u^{p_{[\delta]}}, \alpha_2\in u^p$ and $\alpha_3\in u^{p_{[\delta]}}-u^q$ then $\{\alpha_0,\alpha_1\}\subseteq u^q \wedge \{\alpha_0,\alpha_1\}\nsubseteq u^{p_{[\delta]}}$ while $\{\alpha_1,\alpha_2\}\subseteq u^{p_{[\delta]}} \wedge \{\alpha_1,\alpha_2\}\nsubseteq u^q$.
However, for this to happen we must have a non-monotonic sequence of $\alpha$s, as we must fall back in $u^p$ in the middle of the process.
Hence if $v\subseteq\{0,\ldots,m\}$ and $\langle\alpha_\ell:\ell\in v\rangle$ is increasing then necessarily $\{\alpha_\ell:\ell\in v\}\subseteq u^{p_{[\delta]}}$ or $\{\alpha_\ell:\ell\in v\}\subseteq u^q$.
We remain, therefore, with two cases: \newline

\par\noindent\emph{Case 2}: $\{\alpha_\ell:\ell\in v\}\subseteq u^{p_{[\delta]}}$.

We wish to use the fact that $u^{p_{[\delta]}}$ is a condition and apply it to $v$.
But the fact that $\{\alpha_\ell:\ell\in v\}\subseteq u^{p_{[\delta]}}$ does not guarantee that $\alpha_\ell\in u^{p_{[\delta]}}$ for every $\ell\leq m$, so the assumption of the case is insufficient for this plain argument.
Still, we can argue as follows.
For each $\ell\leq m$ let $\beta_\ell=\alpha_\ell$ if $\alpha_\ell\in u^{p_{[\delta]}}$ and let $\beta_\ell=h(\alpha_\ell)$ if $\alpha_\ell\notin u^{p_{[\delta]}}$ (in which case $\alpha_\ell\in u^q-u^p$).
Notice that $\{\beta_\ell:\ell\leq m\}\subseteq u^{p_{[\delta]}}$.
We claim that $\ell<m\Rightarrow\beta_\ell\neq\beta_{\ell+1}$.

To see this, fix $\ell<m$.
If $\alpha_\ell,\alpha_{\ell+1}\in u^{p_{[\delta]}}$ then $\beta_\ell=\alpha_\ell\neq\alpha_{\ell+1}=\beta_{\ell+1}$.
If $\alpha_\ell,\alpha_{\ell+1}\in u^q$ then $\beta_\ell=h(\alpha_\ell)\neq h(\alpha_{\ell+1})=\beta_{\ell+1}$ since $\alpha_\ell\neq\alpha_{\ell+1}$ and $h$ is one-to-one.
If $\alpha_\ell\in u^{p_{[\delta]}}$ and $\alpha_{\ell+1}\notin u^{p_{[\delta]}}$ then $\alpha_{\ell+1}\in u^q-u^{p_{[\delta]}}$ and we can assume that $\alpha_\ell\notin u^q$ (otherwise $\alpha_\ell,\alpha_{\ell+1}\in u^q$, a case which was covered before), so $\alpha_\ell\in u^{p_{[\delta]}}-u^q$.
But then $\{\alpha_\ell,\alpha_{\ell+1}\}\nsubseteq u^{p_{[\delta]}} \wedge \{\alpha_\ell,\alpha_{\ell+1}\}\nsubseteq u^q$ so this is covered in Case 1.
If $\alpha_\ell\in u^q$ and $\alpha_{\ell+1}\notin u^q$ we use a similar argument.
So we conclude that $\ell<m\Rightarrow\beta_\ell\neq\beta_{\ell+1}$.

Notice further that $f^r(\beta_\ell,\beta_{\ell+1})=f^r(\alpha_\ell,\alpha_{\ell+1})$ whenever $\ell<m$.
Now since $\{\alpha_\ell:\ell\in v\}\subseteq u^{p_{[\delta]}}$ we have $\beta_\ell=\alpha_\ell$ for every $\ell\in v$.
Hence $\{\beta_\ell:\ell\in v\}$ is increasing and $|v|<\max\{f^r(\beta_\ell,\beta_{\ell+1}):\ell<m\} = \max\{f^r(\alpha_\ell,\alpha_{\ell+1}):\ell<m\}$ as $p_{[\delta]}\in\mathbb{Q}$. \newline

\par\noindent\emph{Case 3}: $\{\alpha_\ell:\ell\in v\}\subseteq u^q$.

This is similar to the previous case. \newline

We conclude, therefore, that $r\in\mathbb{Q}$ which shows that $p_{[\delta]}$ and $q$ are compatible.
Having established this general claim let us prove the chain condition.
Assume towards contradiction that $\{q_\alpha:\alpha\in\omega_1\}$ is an antichain in $\mathbb{Q}$.
For every limit ordinal $\alpha\in\omega_1$ let $p_\alpha=q_\alpha\cap\alpha$.
Using Fodor's lemma we see that for some stationary subset $S\subseteq\omega_1$ which consists of limit ordinals and a fixed condition $p$ we have $\alpha\in S\Rightarrow p_\alpha=p$.
We shrink $S$ further by assuming that all the conditions in $\{q_\alpha:\alpha\in S\}$ are pairwise isomorphic.
It means that if $|u^{q_\alpha}|=n_\alpha$ for every $\alpha\in S$ then $n_\alpha=n$ for some fixed $n\in\omega$ and for all the elements of $S$.
Moreover, the order pattern of the elements of every $u^{q_\alpha}$ is assumed to be the same pattern for all the elements of $S$.

Now choose $\zeta,\eta\in S$ so that $\zeta<\eta$ and $u^{q_\zeta}\subseteq\eta$.
Let $\delta$ be a limit ordinal such that $p\subseteq\delta$.
It follows that $p_{\zeta}^{[\delta]}$ and $q_\eta$ are compatible, where $p_{\zeta}^{[\delta]}$ is the projection of $p_\zeta$ to the countable ordinal $\delta$.
But then $q_\zeta$ and $q_\eta$ are also compatible, a contradiction.

\hfill \qedref{lemccc}

We can derive now our conclusion about path relations at $\omega_1$.
We indicate that the above lemma fails while trying to apply the same forcing notion to higher cardinals.
This might serve as an evidence to the possibility that the polarized relation $\binom{\mu^+}{\mu}\rightarrow\binom{\mu^+\ \omega_1}{\mu\ \mu}$ for a singular cardinal with countable cofinality is consistent and maybe even provable under the assumption that $2^\mu=\mu^+$.

\begin{corollary}
\label{coromega1} $\omega_1\nrightarrow_{\rm path}(\omega)^2_{\omega,<\omega}$.
\end{corollary}

\par\noindent\emph{Proof}. \newline
Let $c:[\omega_1]^2\rightarrow\omega$ be as guaranteed in Theorem \ref{thmomega1}.
We claim that $c$ exemplifies the negative relation $\omega_1\nrightarrow_{\rm path}(\omega)^2_{\omega,<\omega}$.
To see this, assume by way of contradiction that $(\alpha_n:n\in\omega)$ is a counterexample.
Namely, there is $m\in\omega$ so that $|\{c(\alpha_\ell,\alpha_{\ell+1}):\ell\in\omega\}|=m$.
Hence there is a natural number $n\in\omega$ such that $c(\alpha_\ell,\alpha_{\ell+1})<n$ for every $\ell\in\omega$.

By induction on $\ell\in\omega$ we choose $j(\ell)\in\omega$ such that if $k<\ell$ then $j(k)<j(\ell)$ and $\alpha_{j(k)}<\alpha_{j(\ell)}$.
Applying the conclusion of Theorem \ref{thmomega1} to the sequence $(\alpha_\ell:\ell\leq j(n))$ and the set $v=\{j(\ell):\ell\leq n\}$ we see that $|v|=n+1<\max\{c(\alpha_\ell,\alpha_{\ell+1}):\ell\leq j(n)\}<n$, a contradiction.

\hfill \qedref{coromega1}

\newpage

\providecommand{\bysame}{\leavevmode\hbox to3em{\hrulefill}\thinspace}
\providecommand{\MR}{\relax\ifhmode\unskip\space\fi MR }
\providecommand{\MRhref}[2]{%
  \href{http://www.ams.org/mathscinet-getitem?mr=#1}{#2}
}
\providecommand{\href}[2]{#2}

\end{document}